\documentclass[10pt]{amsart}
\usepackage{hyperref}
\usepackage{amssymb,amsmath, amsthm}

\newtheorem{thm}{Theorem}[section]

\newtheorem{prop}{Proposition}[section]
\newtheorem{defn}{Definition}[section]

\theoremstyle{definition}

\newcommand{\w}{\omega}        
\newcommand{\J}{\mathbf{J}}    

\newcommand{\Ad}{\mathrm{Ad}}  
\newcommand{\im}{\mathrm{im}\, }  
\newcommand{\FL}{\mathbb{F}\mathrm{L}}  
\newcommand{\Sl}{\mathbf{S}}   
\newcommand{\g}{\mathfrak{g}}
\newcommand{\h}{\mathfrak{h}}
\newcommand{\be}{\begin{equation}}
\newcommand{\ee}{\end{equation}}
\newcommand{\bea}{\begin{eqnarray}}
\newcommand{\eea}{\end{eqnarray}}
\newcommand{\restr}[1]{\vrule height3ex width.4pt depth1.4ex\lower1.4ex\hbox{\scriptsize $\,#1$}}
\newcommand{\rrestr}[1]{\vrule height2ex width.4pt depth0.9ex\lower0.9ex\hbox{\scriptsize $\,#1$}}

\begin{document}

\title{The isotropy lattice of a lifted action}

\author{Miguel Rodr\'{i}guez-Olmos}
\thanks{M. Rodr\'{i}guez-Olmos: Section de Math\'ematiques, EPFL,
CH-1015 Lausanne, Switzerland. miguel.rodriguez@epfl.ch}


\begin{abstract}
We obtain a characterization of the isotropy lattice for the
lifted action of a Lie group $G$ on $TM$ and $T^*M$ based only on
the knowledge of $G$ and its action on $M$. Some applications to
symplectic geometry are also shown.
\end{abstract}

\maketitle
\section{Introduction}
In this note we make some remarks about the geometry of
(co-)tangent-lifted actions that seem to be unknown or not
available in the literature. The problem considered is the
following: if a Lie group $G$ acts properly on the manifold $M$,
this action is characterized by its isotropy lattice, which is the
set of conjugacy classes of stabilizers of points of $M$,
partially ordered by the relation of subconjugation. The knowledge
of this lattice is an important tool that offers valuable
information about the topology of the stratification  of the
quotient space $M/G$ and its singularities. In many cases of
geometrical and mechanical importance, one is interested in the
study of the quotient space of a tangent (or cotangent) bundle by
the lift of a proper group action on its base. In these cases, the
relevant information is then obtained by the knowledge of the
isotropy lattice for the (co-)tangent-lifted action. Now, since
both the tangent bundle of a manifold $M$ and the tangent-lifted
action of a group $G$ on $TM$ are completely obtained from the
geometry of $M$ and the $G$-action on it, the isotropy lattice for
the lifted action should be obtainable from the isotropy lattice
for the base action (supposed to be known) without the need of a
direct study of the tangent-lifted action separately. Our main
result, Theorem \ref{mainthm}, gives a characterization of the
isotropy lattice for the lifted action of $G$ on $TM$ once the
isotropy lattice for the action on $M$ and the adjoint
representation of $G$ are known. From this result, an algorithm to
obtain every element in the lattice of the lifted action can be
easily implemented. Also, we obtain several conditions relating
stabilizers and orbit types of $TM$ and $M$, and some applications
to the study of level sets of the momentum map in symplectic
geometry.

Throughout this note $M$ will denote a smooth, paracompact,
connected and finite-dimensional manifold on which the
finite-dimensional Lie group $G$ acts smoothly and properly. This
action induces in the usual way a proper action on the tangent
bundle $TM$ (respectively on $T^*M$) which is also smooth and
proper.
\begin{defn}\label{defproper} If
$H\subseteq G$ is a subgroup of $G$ and $x\in M$, denote by $(H)$
the conjugacy class of $H$ in $G$ and by $G_x$ the stabilizer of
$x$.
\begin{itemize}
\item[(i)] The set $I(G,M)=\{(H)\,:\, \mathrm{there\,\, is\,\,
some}\,\, x\in M$ with $G_x=H\}$ is called the {\bf isotropy
lattice} for the $G$-action on $M$. \item[(ii)] For any $(H)\in
I(G,M)$, the set $M_{(H)}=\{x\in M\,:\, G_x\in (H)\}$ is called an
{\bf orbit type} (of type $H$) of $M$.
\end{itemize}
\end{defn}
Sometimes we will also consider conjugacy classes with respect to
smaller subgroups of $G$. We will use the same symbol $(H)$ for
these new classes without explicitly saying it, unless this is not
clear from the context. Also, it is clear from the definition the
meaning of $I(G,X)$ and $X_{(H)}$ for any subspace $X$ of $M$.

 There is a partial order on an isotropy lattice as follows:
if $(H_1)$ and $(H_2)$ are two different elements of  $I(G,M)$
then $(H_1)< (H_2)$ if $H_1$ is conjugated to a proper subgroup of
$H_2$. This condition is easily checked to be independent of the
representatives chosen.
 The following proposition collects a number
of important properties of this action, which can be consulted for
example in \cite{DuiKol,OrRa2003,Palais}.
\begin{prop}\label{prelimprop}For a $G$-action on $M$ the following hold:
\begin{itemize}
\item[(i)] Every stabilizer group $G_x$ is compact.
\item[(ii)]There is a $G$-invariant metric on $M$, \item[(iii)]
(Tube Theorem) Let $\Sl$ be a $G_x$-invariant orthogonal
complement to $\g\cdot x$ in $T_xM$ with respect to some
$G$-invariant Riemannian metric on $M$ and $O\subset \Sl$ a small
enough $G_x$-invariant open ball centered at the origin in $\Sl$.
Then  the space $G\times_{G_x}O$ obtained by quotienting $G\times
O$ by the $G_x$-action $g'\cdot (g,s)=(g{g'}^{-1},g'\cdot s)$ is
$G$-diffeomorphic to a $G$-invariant open neighborhood $U$ of the
orbit $G\cdot x$ by the map $\phi(\left[g,s\right])= g\cdot \exp_x
(s),$ where the $G$-action on the left hand side is given by
$g_1\cdot\left[g_2,s\right]=\left[g_1g_2,s\right]$. \item[(iv)]
The connected components of each orbit type $M_{(H)}$ are embedded
submanifolds of $M$, and $M_{(H)}\subset \overline{M_{(H')}}$ for
every $(H')\leq (H)$. \item[(v)] There is a minimal class $(H_0)$
in $I(G,M)$ such that $(H_0)\leq (H)$ for every $(H)\in I(G,M)$.
\end{itemize}
\end{prop}

{\flushleft {\bf Remark 1.}}\label{remcoadjoint} Besides $I(G,M)$,
we will consider the isotropy lattice for the $H$-representation
on $\g/\h$ given by $h\cdot [\xi]=\left[\Ad_h\xi\right]$ for a
compact subgroup $H$. Note that this action is induced from the
restriction to $H$ of the adjoint representation of $G$ on $\g$,
and that the dual of $\g/\h$ is isomorphic to the annihilator of
$\h$, $\mathrm{ann}\, \h\subset \g^*$. Furthermore, since $H$ is
compact this isomorphism can be chosen $H$-equivariant with
respect to the restricted coadjoint representation
$H\times\mathrm{ann}\,\h\rightarrow \mathrm{ann}\,\h$ given by
$h\cdot\mu=(\Ad_{h}^*)^{-1}\mu$. Therefore, pairs of elements
identified by this isomorphism have identical stabilizers and
hence $I(H,\g/\h)=I(H,\mathrm{ann}\,\h)$.

\section{The Main Result}
\begin{thm}\label{mainthm}
Let $G$ act on $M$ and by tangent lifts on $TM$, and $L$ be a
subgroup of $G$. Then $(L)\in I(G,TM)$ if and only if there exist
$(H_1),(H_2)\in I(G,M)$ and $(K)\in I(H_2,\mathrm{ann}\,\h_2)$
such that $(H_1)\leq (H_2)$ and $L=H_1\cap K$.
\end{thm}
{\bf Proof}. We will fix once and for all a $G$-invariant
Riemannian metric on $M$. Let $x\in M$ have stabilizer $H=G_x$,
then we can form the $H$-invariant orthogonal splitting
$T_xM=\g\cdot x\oplus\Sl$ as in the Tube Theorem (item (iii) in
Proposition \ref{prelimprop}). Note that $\xi_M(x)=0$ if and only
if $\xi\in\h$, and so there is an $H$-isomorphism $\psi
:\g/\h\times\Sl \rightarrow T_xM$ given by $\psi ([\xi],s)=
\xi_M(x)+s$. $H$-equivariance is understood with respect to the
induced linear action of $H$ on $T_xM$ and to the diagonal action
on $\g/\h\times\Sl$ given by $h\cdot
([\xi],s)=(\left[\Ad_h\xi\right],h\cdot s)$. \

By equivariance of the tangent bundle projection
$\tau:TM\rightarrow M$, for any element $v_x\in TM$ such that
$\tau(v_x)=x$ the tangent bundle and base isotropy groups are
related by $G_{v_x}\subseteq H=G_x$. Furthermore,
$G_{v_x}=H_{v_x}$, where on the left-hand side we consider the
stabilizer for the full lifted action of $G$ on $TM$ and on the
right-hand side the stabilizer for the induced linear
representation of $H$ on the vector space $T_xM$. If
$\psi^{-1}(v_x)=([\xi],s)$ then $G_{v_x}=H_{[\xi]}\cap H_s$, which
means that the stabilizer of any tangent vector $v_x$ over $x$ is
the intersection of two representatives of the lattices
$I(H,\g/\h)$ and $I(H,\Sl)$. Conversely, any such intersection is
the stabilizer of some vector in $T_xM$. Recall now from Remark
\ref{remcoadjoint} that $I(H,\g/\h)=I(H,\mathrm{ann}\,\h)$. We
need now to identify the elements of $I(H,\Sl)$. For that, using
the Tube Theorem choose elements $s\in O\subset\Sl$ and $g\in G$.
Then
\be\label{conjugisotropy}gH_sg^{-1}=G_{[g,s]}=G_{\phi([g,s])}.\ee
Note that using every $g\in G$ and $s\in O$ the image of $\phi$
covers the full neighborhood $U\subset M$ in the Tube Theorem.
Also, $x\in M_{(H)}$, and then by item (iv) in Proposition
\ref{prelimprop} we have that for any $(H')\in I(G,M)$, $U\cap
M_{(H')}\neq\emptyset\Leftrightarrow (H')\leq (H)$. This, together
with \eqref{conjugisotropy}, means that the stabilizers for the
linear $H$-representation on $\Sl$ are conjugated to subgroups
$H'$ such that $(H')\in I(G,M)$ and $(H')\leq (H)$. Conversely, if
$(H')\in I(G,M)$ satisfies $(H')\leq (H)$, then there is a
representative $H'$ of $(H')$ with $H'\subset H$ and such that the
$H$-conjugacy class $(H')$ belongs to $I(H,\Sl)$. Now, making
$H=H_2$ and $H'=H_1$ the result follows. \hfill $\Box$

{\flushleft {\bf Remark 2.}} We can obtain all the classes in
$I(G,TM)$ with an easy algorithm. First, choose representatives of
every class in $I(G,M)$ such that if $(H')\leq (H)$ the
corresponding representatives satisfy $H'\subset H$. It is always
possible to choose a complete set of representatives for all the
classes of $I(G,M)$ in this way, and we will call them normal
representatives. Let $H_0$ be the normal representative
corresponding to the minimal class of $I(G,M)$ (see item (v) in
Proposition \ref{prelimprop}). We will say that $H_0$ has depth
zero. Any other normal representative $H$ will have depth $n+1$ if
there is a normal representative $H'$ with depth $n$ such that 1)
$H'\subset H$ and 2) there is no other normal
representative $H''$ with $H'\subset H''\subset H$.\\
 To compute
all the classes in $I(G,TM)$ we start by computing the classes of
$I(H_0,\mathrm{ann}\,\h_0)$. Then, for any $n$ and every normal
representative $H$ of depth $n$ intersect the classes of
$I(H,\mathrm{ann}\,\h)$ with the $H$-class of every normal
representative of depth $0,\ldots,n-1$ included in $H$. All the
classes obtained after iterating in $n$ in this way can be made
$G$-classes by con\-ju\-ga\-ting in $G$. After removing the
repeated ones, we obtain all the elements of $I(G,TM)$.

{\flushleft {\bf Remark 3.}}\label{alsocotangent}
\begin{itemize}\item[(i)] Note that the choice of a $G$-invariant metric on $M$
(always available by item (ii) in Proposition \ref{prelimprop})
 provides a $G$-bundle isomorphism $\FL:TM\rightarrow T^*M$.
  Consequently, $I(G,TM)=I(G,T^*M)$, and since $\FL$
covers the identity in $M$, $\tau \left( (TM)_{(H)}\right)=\tau
\left( (T^*M)_{(H)}\right)$ for any $(H)\in I(G,TM)$. Here we have
denoted also by $\tau$ the cotangent bundle projection
$T^*M\rightarrow M$. \item[(ii)] It easily follows from the proof
of Theorem \ref{mainthm} that if $G$ is Abelian, then
$I(G,TM)=I(G,M)$. Furthermore, for any $(H)\in I(G,M)$, then
$\tau\left((TM)_{(H)}\right)=\overline{M_{(H)}}$.
\end{itemize}

\section{Applications in symplectic geometry}
Using Remark \ref{alsocotangent} we can translate every result
previously obtained for $I(G,TM)$ to $I(G,T^*M)$ where $G$ acts on
$T^*M$ by cotangent lifts. Recall that this lifted action is
Hamiltonian for the canonical symplectic structure $\w$ on $T^*M$,
and then it is natural to study the properties of $I(G,T^*M)$ from
the point of view of symplectic geometry. In this section we make
some remarks along these lines. Namely, we compute the restricted
isotropy lattice for every $G$-invariant level set of the momentum
map. The knowledge of the isotropy lattice of level sets of the
momentum map is an important piece of information in the theory of
singular reduction (see \cite{OrRa2003,SaLe}). From the dynamics
side we also obtain a restriction on the conjugacy classes of
stabilizers for possible relative equilibria in an important class
of Hamiltonian systems.

Recall that a cotangent-lifted action of $G$ on $(T^*M,\w$) is
Hamiltonian with equivariant momentum map $\J:T^*M\rightarrow\g^*$
defined by
$\langle\J(\alpha_x),\xi\rangle=\langle\alpha_x,\xi_M(x)\rangle$,
for every $\alpha_x\in T^*_xM$ and $\xi\in\g$. This means that a
necessary condition for $\alpha_x\in \J^{-1}(\mu)$ is
$\mu\in\text{ann}\,\g_x$. Let now $\mu\in\text{im}\,\J\subset\g^*$
be a totally isotropic momentum value, i.e. $G_\mu=G$ for the
coadjoint representation. Equivariance of $\J$ implies that
$\J^{-1}(\mu)$ is $G$-invariant. Note also that since $\mu$ is
totally isotropic, if $\mu\in\text{ann}\,\g_x$ the same is true
for the Lie algebra of any other representative of $(G_x)$. We
will define the $\mu$-lattice of $M$ as
$$I^\mu (G,M)=\{(H)\in I(G,M)\,:\,\mu\in\text{ann}\,\h\}.$$ If
$(H)\in I^\mu (G,M)$ then the $\mu$-closure of $M_{(H)}$ is
defined as
$$\text{cl}^\mu\,(M_{(H)})=\left\{\coprod_{(K)}M_{(K)}\,:\, (K)\in
I^\mu(G,M)\,\,\text{and}\, (K)\geq (H)\right\}.$$ Note that any
nonempty $\mu$-closure is $G$-invariant.

Recall also that a Hamiltonian system on $(T^*M,\w)$ acted by $G$
by cotangent lifts is called a (symmetric) ``simple mechanical
system'' if its Hamiltonian function is of the form
$h(\alpha_x)=\frac 12\|\FL^{-1} (\alpha_x)\|^2+V(x)$, where $V$ is
a $G$-invariant function on $M$ and the norm $\|\cdot\|$ and the
Legendre map $\FL$ (see Remark \ref{alsocotangent}) are taken with
respect to a given Riemannian metric for which the $G$-action on
$M$ is isometric. It is well known (se \cite{MarLec}) that if
$\alpha_x$ is a relative equilibrium\footnote{A relative
equilibrium in a Hamiltonian system is a point for which its
Hamiltonian evolution lies inside a group orbit.} for this
Hamiltonian system, then it must be of the form
$\alpha_x=\FL(\xi_M(x))$ for some $\xi\in\g$. We will then call
the $G$-invariant subset of $T^*M$ defined as
$\{\FL(\xi_M(x))\,:\,x\in M\,\text{and}\,\xi\in\g\}$ the set of
``possible relative equilibria''. Note that this set depends only
on $M$ and the $G$-action on it, and is the same for any simple
mechanical system defined on $T^*M$. The knowledge of the
restricted isotropy lattice for the set of possible equilibria on
$T^*M$ gives an estimate of the stabilizers that relative
equilibria of simple mechanical systems on $T^*M$ can have. This
fact is of importance in the qualitative analysis of these
systems, as for instance it can predict the existence or
exclusions of certain types of equivariant bifurcations for any
simple mechanical system defined on $T^*M$.

\begin{prop}
For the cotangent lift of $G$ on $T^*M$,  and if
$\mu\in\mathfrak{g}^*$ is any totally isotropic momentum value,
the following properties are satisfied:
\begin{itemize} \item[(i)] $I(G,\J^{-1}(0))=I(G,M)$ and for any
$(H)\in I(G,M)$, $\tau\left( \left(  \J^{-1}(0)
\right)_{(H)}\right)=\overline{M_{(H)}}$.\item[(ii)]$I(G,\J^{-1}(\mu))=I^\mu(G,M)$
and for any $(H)\in I^\mu(G,M)$, $\tau\left( \left(  \J^{-1}(\mu)
\right)_{(H)}\right)=\mathrm{cl}^\mu\left(M_{(H)}\right)$.\item[(iii)]
For any $G$-invariant simple mechanical system on $M$, the
isotropy lattice of the set of possible relative equilibria is
$$\coprod_{(H)\in I(G,M)}
G\cdot I(H,\mathrm{ann}\,\h).$$ Where $G\cdot
I(H,\mathrm{ann}\,\h)$ is the set of conjugacy classes in $G$ of
representatives of elements of $I(H,\mathrm{ann}\,\h)$.
\end{itemize}
\end{prop}
{\bf Proof} $(i)$ is clearly a particular case of $(ii)$. To prove
$(ii)$ fix $x\in M$ with $G_x=H$. Recall from the definitions of
$\J,\,\FL$ and $\psi:\g/\h\times\Sl\rightarrow T_xM$ that if
$\alpha_x\in\J^{-1}(\mu)\cap T^*_xM$ then
$\alpha_x=(\FL\circ\psi)([\xi],s)$, where $[\xi]$ is fixed by $H$.
This implies that $G_{\alpha_x}=H_s$, which already appears as an
stabilizer of some point of $M\cap \im (\exp_x)$. Also since
$H_s\subseteq H$ and $\mu\in\text{ann}\,\h$, then
$\mu\in\text{ann}\,\g_{\alpha_x}$ as well. Applying this to every
point $x\in M$ we find that $I(G,\J^{-1}(\mu))=I^\mu(G,M)$. Also,
the same reasoning shows that for each $(K)\in I^\mu(G,M)$ such
that $(K)\geq (H)$, there is at least an element in
$(\J^{-1}(\mu))_{(H)}\cap T^*_xM$ at every $x\in M_{(K)}$, which
ends the proof of $(ii)$.

For $(iii)$, recall that by equivariance of the map $\FL$, the
isotropy lattice for the set of possible relative equilibria is
the same as the isotropy lattice for the subset of $TM$ given by
the collection of infinitesimal generators for the $G$-action on
$M$. Fixing again a point $x\in M$ with isotropy $H=G_x$, every
infinitesimal generator at $x$ is of the form $v_x=\xi_M(x)$ or
equivalently $\psi([\xi],0)$ and hence the $H$-class $(G_{v_x})$
is in $I(H,\text{ann}\,\h)$, by the proof of Theorem
\ref{mainthm}. Computing this lattice over each base point in
$M_{(H)}$ generates $G\cdot I(H,\text{ann}\,\h)$, since any such
base point has stabilizer conjugated to $H$. Finally, doing the
same over each orbit type of $M$ is equivalent to taking the
reunion of all $G\cdot I(H,\text{ann}\,\h)$ for every $(H)\in
I(G,M)$.\hfill $\Box$



\begin{thebibliography}{00}

\bibitem{DuiKol} J.J. Duistermaat and J.A.C. Kolk [2000], {\it Lie
groups}, Universitext, Springer-Verlag.


\bibitem{MarLec} J.E. Marsden [1992], {\it Lectures on Mechanics}, Lecture Note Series
{\bf 174}, LMS, Cambridge University Press.

\bibitem{OrRa2003} J.-P. Ortega and T.S. Ratiu [2004], {\it Momentum maps and Hamiltonian reduction},
Progress in Mathematics {\bf 222}, Birkhauser-Verlag.

\bibitem{Palais} R.S. Palais [1961], On the existence of slices for
non-compact Lie group, {\it Ann. of Math.}, {\bf 73}, 265--323.

\bibitem{SaLe} R. Sjamaar and E. Lerman [1991], Stratified symplectic
spaces and re\-duc\-tion, {\it Ann. of Math.} {\bf 134}, 375--422.


\end{thebibliography}
 \end{document}